\documentclass[a4paper]{llncs}

\usepackage[utf8x]{inputenc}

\usepackage{amsmath}

\usepackage{amsthm}
\usepackage{amssymb}
\usepackage{breqn}
\usepackage{framed}
\usepackage{array}
\usepackage{xfrac}
\usepackage{enumerate}
\usepackage{bookmark}
\usepackage{hyperref}
\usepackage{graphicx}

\hypersetup{
  pdftitle = {Digits of pi: limits to the seeming randomness},
  pdfauthor = {Karlis Podnieks},
  pdfkeywords = {digits of pi, random digits},
  pdfstartview=FitH,
  unicode=true
}

\begin{document}

\bookmarksetup{startatroot}

\title{Digits of pi: limits to the seeming randomness}
\author{Karlis Podnieks}

\institute{University of Latvia, Raina bulvaris 19, Riga, LV-1586, Latvia}

\theoremstyle{plain}
\newtheorem{hypo}{Hypothesis}
\newtheorem{cor}{Corollary}
\newtheorem{que}{Question}
\newtheorem{obs}{Observation}

\maketitle

\bookmarksetup{startatroot}

\begin{abstract}

The decimal digits of $\pi$ are widely believed to behave like as statistically independent random variables taking the values $0, 1, 2, 3, 4, 5$, $6, 7, 8, 9$ with equal probabilities $1/10$.

In this article, first, another similar conjecture is explored - the seemingly almost random behaviour of digits in the base 3 representations of powers $2^n$. This conjecture seems to confirm well - it passes even the tests inspired by the Central Limit Theorem and the Law of the Iterated Logarithm.

After this, a similar testing of the sequences of digits in the decimal representations of the numbers $\pi$, $e$ and $\sqrt{2}$ was performed. The result looks surprising: unlike the digits in the base 3 representations of $2^n$, instead of oscillations with amplitudes required by the Law of the Iterated Logarithm, {\itshape convergence to zero} is observed. If, for such "analytically" defined irrational numbers, the observed behaviour remains intact {\itshape ad infinitum}, then the seeming randomness of their digits is only a limited one.

\end{abstract}

\pdfbookmark[1]{1. Introduction}{intro}
\section{Introduction}
\label{intro}
The decimal digits of $\pi$ are widely believed to behave like as statistically independent random variables taking the values $0, 1, 2, 3, 4, 5, 6, 7, 8, 9$ with equal probabilities $\frac{1}{10}$ (for an overview, see \cite{marsaglia}). 

In Section \ref{powers} (it reproduces in part Section 2.1 of \cite{paper2}) another similar conjecture is explored - the seemingly almost random behaviour of digits in the base 3 representations of powers $2^n$. This conjecture seems to  confirm well - it passes even the tests inspired by the Central Limit Theorem and the Law of the Iterated Logarithm.
 
Especially remarkable is Fig.\ref{Oscilat1} below showing oscillations with amplitudes almost as required by the Law of the Iterated Logarithm.

In Section \ref{pi}, similar pictures for the sequences of digits in the decimal representations of the numbers $\pi$, $e$ and $\sqrt{2}$ are obtained. The result looks surprising: instead of oscillations with amplitudes required by the Law of the Iterated Logarithm, we observe {\itshape convergence to zero}!  If, for such "analytically" defined irrational numbers, the observed behaviour remains intact {\itshape ad infinitum}, then the seeming randomness of their digits is only a limited one.

\pdfbookmark[1]{2. Base 3 representations of powers of 2}{powers}
\section{Base 3 representations of powers of 2}
\label{powers}

Throughout this section, it is assumed that $p, q$ are positive integers such that  $\frac{\log{p}} {\log{q}}$ is irrational, i.e., $p^a \neq q^b$ for any integers $a, b>0$.

\begin{definition}
Let us denote by $D_q(n, i)$ the \textbf{$i$-th digit} in the canonical base $q$ representation of the number $n$, and by $S_q(n)$ - the \textbf{sum of digits} in this representation. 
\end{definition}

Let us consider base $q$ representations of powers $p^n$. Imagine, for a moment, that, for fixed $p, q, n$,  most of the digits  $D_q(p^n, i)$ behave like as statistically independent random variables taking the values $0, 1, ..., q-1$ with equal probabilities $\frac{1}{q}$. Then, the (pseudo) mean value and (pseudo) variance of $D_q(p^n, i)$ should be
 \[E=\frac{q-1}{2}; V=\sum\limits_{i=0}^{q-1}\frac{1}{q}\left(i-\frac{q-1}{2}\right)^2=\frac{q^2-1}{12}.\]
The total number of digits in the base $q$ representation of $p^n$ is $k_n \approx n \log_q{p}$, hence, the (pseudo) mean value of the sum of digits  $S_q(p^n)=\sum \limits_{i=1}^{k_n} D_q(p^n, i)$ should be $E_n \approx n\frac{q-1}{2}\log_q{p}$ and, because of the assumed (pseudo) independence of digits, its (pseudo) variance should be $V_n \approx n\frac{q^2-1}{12}\log_q{p}$. As the final consequence, the corresponding centered and normed variable
\[\frac{S_q(p^n)-E_n}{\sqrt{V_n}}\]
should tend to behave as a standard normally distributed random variable with probability density $\frac{1}{\sqrt{2\pi}}e^{-\frac{x^2}{2}}$.

\begin{figure}
    \centering
    \includegraphics[width=0.8\textwidth]{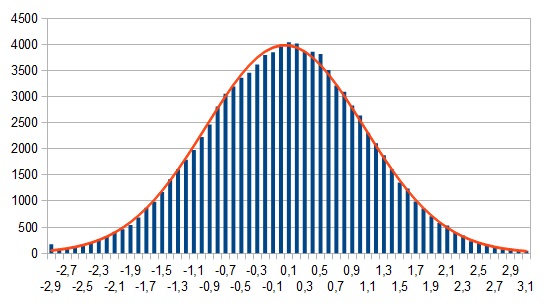}
    \caption{Histogram of the centered and normed variable $s_3(2^n)$}
    \label{Histo1}
\end{figure}

One can try to verify this conclusion experimentally. For example, let us compute $S_3(2^n)$ for $n$ up to $100000$, and let us draw the histogram of the corresponding centered and normed variable
\[s_3(2^n)=\frac{S_3(2^n)-n \log_3 2}{\sqrt{n\frac{2}{3}\log_3 2 }}\]
(see Fig. \ref{Histo1}). As we see, this variable behaves, indeed, almost exactly, as a standard normally distributed random variable (the solid curve).

Observing such a phenomenon ``out there'', one could conjecture that the sum of digits $S_q(p^n)$, as a function of $n$, behaves almost as $n\frac{q-1}{2}\log_q{p}$, i.e., almost \textbf{linearly} in $n$.

An even more advanced idea for testing randomness of sequences of digits was proposed in \cite{borwein1} - let us use the Law of the Iterated Logarithm. Namely, let us try to estimate the amplitude of the possible deviations of  $S_q(p^n)$ from $n\frac{q-1}{2}\log_q{p}$ by ``applying'' the Law of the Iterated Logarithm.

Let us consider the following centered and normed (pseudo) random variables:
\[d_q(p^n, i) = \frac{D_q(p^n, i)-\frac{q-1}{2}}{\sqrt{\frac{q^2-1}{12}}}.\]
By summing up these variables for $i$ from $1$ to $k_n$, we obtain a sequence of (pseudo) random variables:
\[\kappa_q(p, n) = \frac{S_q(p^n)-\frac{q-1}{2}k_n}{\sqrt{\frac{q^2-1}{12}}},\]
that ``must obey" the Law of the Iterated Logarithm. Namely, if the sequence $S_q(p^n)$ behaves, indeed, as a "typical" sum of equally distributed random variables, then $\lim \limits_{n \to \infty} \inf$ and $\lim \limits_{n \to \infty} \sup$ of the fraction
\[\frac{\kappa_q(p, n)}{\sqrt{2 k_n \log{\log k_n}}},\]
 must be $-1$ and $+1$ correspondingly ($\log$ stands for the natural logarithm).

Therefore, it seems, we could conjecture that, if
\[\delta_q(p, n) = \frac{S_q(p^n)-(\frac{q-1}{2}\log_q p)n}{\sqrt{(\frac{q^2-1}{6}\log_q p) n \log{\log n}}},\]
then
\[\lim \limits_{n \to \infty} \sup \delta_q(p, n)=1; \lim \limits_{n \to \infty} \inf \delta_q(p, n) = -1.\]

In particular, this would mean that
\[S_q(p^n)=(\frac{q-1}{2}\log_q p)n +O(\sqrt{n \log \log n}).\]
And, for $p=2; q=3$ this would mean (note that $\log_3 2 \approx 0.6309$):
\[S_3(2^n)=n \cdot \log_3 2 + O(\sqrt{n \log \log n});\]
\[\delta_3(2, n) =\frac{S_3(2^n)-n \log_3 2}{\sqrt{(\frac{4}{3} \log_3 2) n \log{\log n}}}\approx \frac{S_3(2^n)-0.6309 n}{\sqrt{0.8412 n  \log \log n}},\] 
\[\lim \limits_{n \to \infty} \sup \delta_3(2, n)=1; \lim \limits_{n \to \infty} \inf \delta_3(2, n) = -1.\]
 
\begin{figure}
    \centering
    \includegraphics[width=0.8\textwidth]{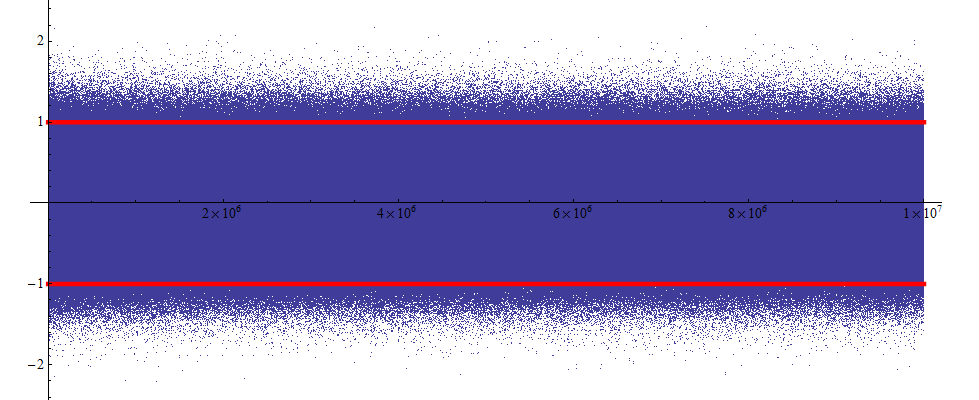}
    \caption{Oscillating behaviour of the expression $\delta_3(2, n)$}
    \label{Oscilat1}
\end{figure}

However, the real behaviour of the expression $\delta_3(2, n)$ until $n=10^7$ does not show convergence of oscillations to the segment $[-1, +1]$ (see Fig. \ref{Oscilat1}, obtained by Juris Čerņenoks). Although $\delta_3(2, n)$ is oscillating almost as required by the Law of the Iterated Logarithm, very many of its values lay outside the segment $[-1, 1]$.

Could we hope to prove the above estimate of $S_q(p^n)$? To my knowledge, the best result on this problem is due to C. L. Stewart. It follows from his Theorem 2 in  \cite{stewart} (put $\alpha=0$), that
\[S_q(p^n)>\frac{\log{n}}{\log{\log{n}} +C_0}-1,\]
where the constant $C_0>0$ can be effectively computed from $q, p$. Since then, no better than $\frac{\log{n}}{\log\log{n}}$ lower bounds of $S_q(p^n)$ have been  proved.

\pdfbookmark[1]{3. Digits of $pi$, $e$ and $sqrt2$}{pi}
\section{Digits of $\pi$, $e$ and $\sqrt{2}$}
\label{pi}

In Section \ref{powers} above, the Central Limit Theorem (Fig. \ref{Histo1}) and the Law of the Iterated Logarithm (Fig. \ref{Oscilat1}) were used to verify the conjecture that the sum of digits of the base 3 representation of $2^n$ behaves closely to the expected behaviour of the sum of the first $n$ members of a sequence of independent random variables taking the values $0, 1, 2$ with equal probabilities $\frac{1}{3}$.

Let us try, as proposed in \cite{borwein1}, to apply this method  to the sequences of digits in the decimal representations of the numbers $\pi$, $e$ and $\sqrt{2}$.
 
Imagine, for a moment, that the digits in the decimal representation of some real number X behave, indeed, like as statistically independent random variables taking the values $0, 1, 2, 3, 4, 5, 6, 7, 8, 9$ with equal probabilities $\frac{1}{10}$. Then, the (pseudo) mean value and (pseudo) variance of $n$-th digit would be (see the formulas above)
 $\frac{10-1}{2}=4.5$ and $\frac{10^2-1}{12} = 8.25$ correspondingly.
And, the (pseudo) mean value of the sum of the first $n$ digits  $S_{10}(n)$ would be $4.5 n$, and, because of the assumed (pseudo) independence of digits, its (pseudo) variance would be $8.25 n$. Let us try to estimate the amplitude of the possible deviations of $S_{10}(n)$ from the expected mean $4.5n$ by ``applying" the Law of the Iterated Logarithm. Let us introduce the necessary centered and normed (pseudo) random variables:
\[\frac{d(i)-4.5}{\sqrt{8.25}}\]
($d(i)$ denotes the $i$-th digit). By summing up these variables for $i$ from $1$ to $n$, we obtain a sequence of (pseudo) random variables:
\[\frac{S_{10}(n)-4.5 n}{\sqrt{8.25}},\]
that ``must obey" the Law of the Iterated Logarithm. Namely, if the sequence $S_{10}(n)$ behaves, indeed, as a "typical" sum of equally distributed random variables, then $\lim \limits_{n \to \infty} \inf$ and $\lim \limits_{n \to \infty} \sup$ of the fraction
\[\delta(n) = \frac{S_{10}(n)-4.5 n}{\sqrt{2\cdot 8.25 n \log{\log n}}},\]
  must be $-1$ and $+1$ correspondingly ($\log$ stands for the natural logarithm).

In particular, this would mean that
\[S_{10}(n)=4.5 n +O(\sqrt{n \log \log n}),\]
and that the values of $\delta(n)$ must oscillate accross the entire segment $[-1, +1]$, like as in Fig. \ref{Oscilat1}. 

\begin{figure}
    \centering
    \includegraphics[width=0.8\textwidth]{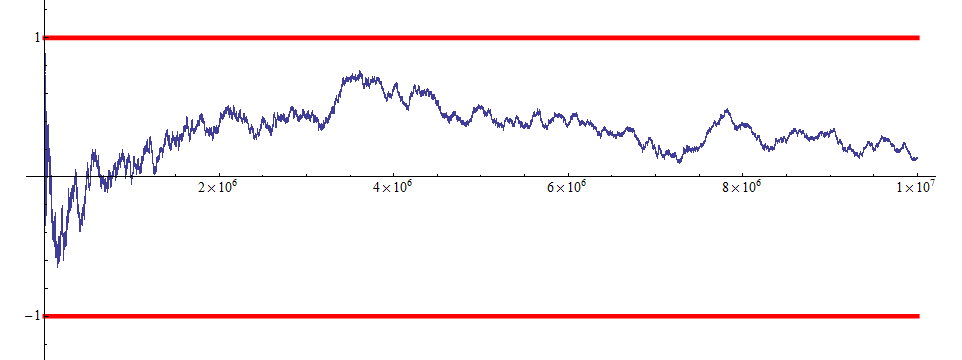}
    \caption{The number $\pi$: behaviour of the expression $\delta(n)$}
    \label{OscilatPI}
\end{figure}

\begin{figure}
    \centering
    \includegraphics[width=0.8\textwidth]{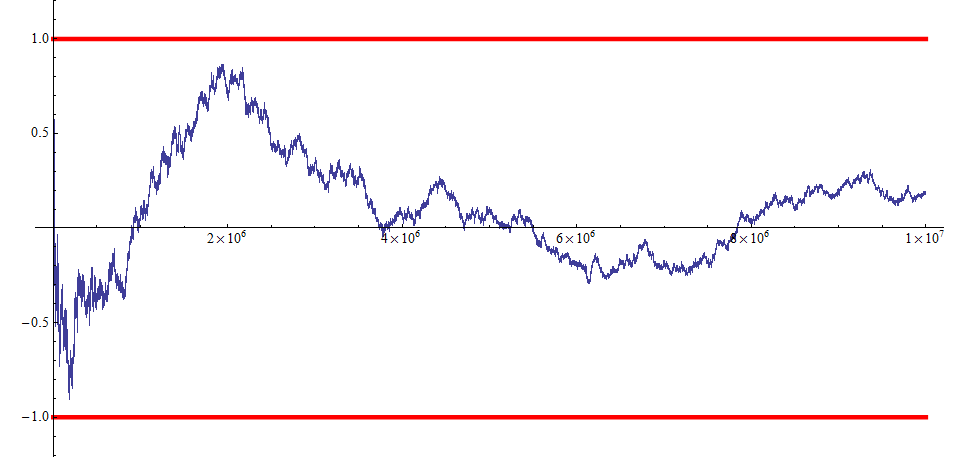}
    \caption{The number $e$: behaviour of the expression $\delta(n)$}
    \label{OscilatE}
\end{figure}

\begin{figure}
    \centering
    \includegraphics[width=0.8\textwidth]{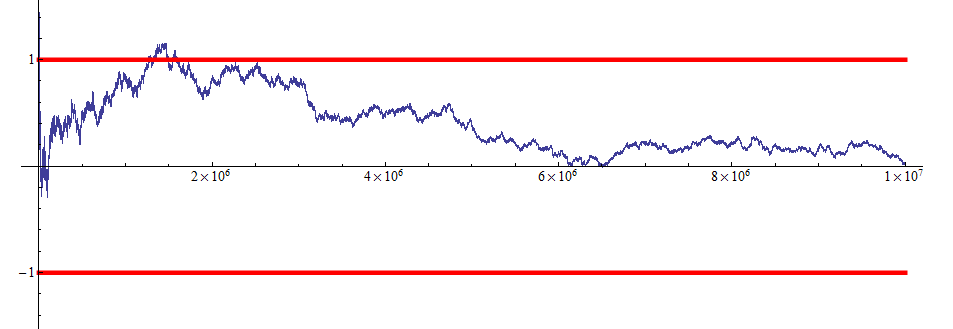}
    \caption{The number $\sqrt{2}$: behaviour of the expression $\delta(n)$}
    \label{OscilatSQRT2}
\end{figure}

However, Fig. \ref{OscilatPI}, Fig. \ref{OscilatE}, and Fig. \ref{OscilatSQRT2} obtained by Juris Čerņenoks for the first $10^7$ digits of the numbers $\pi$, $e$ and $\sqrt{2}$ are showing a completely different behaviour. (Digits were provided by Wolfram Mathematica \cite{wolfram}.)

For the three numbers in question, the values of $\delta(n)$ do not oscillate accross the entire segment $[-1, +1]$, instead, they seem converging to 0. Thus, the pictures seem to support the following somewhat stronger conjecture for $\pi$, $e$ and $\sqrt{2}$ :
\[S_{10}(n)=4.5 n +o(\sqrt{n \log \log n}).\]

\begin{figure}
    \centering
    \includegraphics[width=0.8\textwidth]{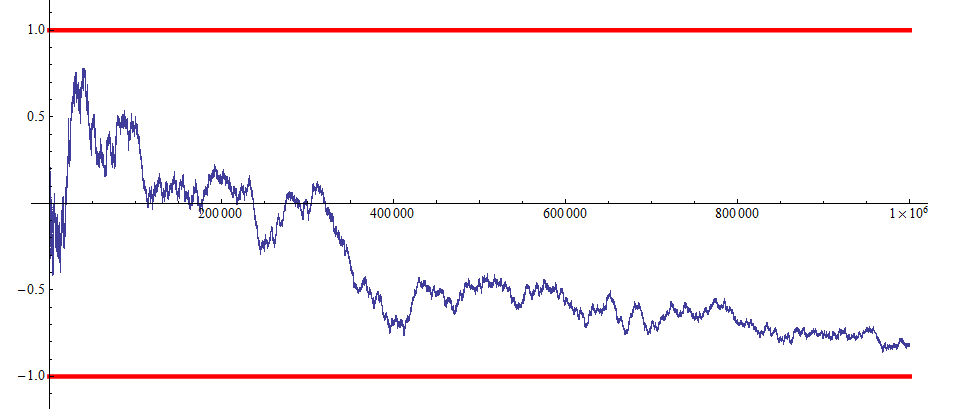}
    \caption{Million random digits from RAND Corp.: behaviour of the expression $\delta(n)$}
    \label{OscilatRAND}
\end{figure}

An even more specific behaviour are showing (see Fig. \ref{OscilatRAND}) the famous Million Random Digits of the RAND Corporation published in 1955 \cite{rand}.

The pictures obtained for $\pi$, $e$, $\sqrt{2}$ are similar to Figure 12(b) obtained for the number $\alpha_{2,3}$ by the authors of \cite{borwein2}. They conclude:

``For  $\alpha_{2,3}$, the corresponding computation of the first $10^9$ values of  $\frac{m_1(n)-n/2}{\sqrt{2n\log\log n}}$ leads to the plot in Figure 12(b) and leads us to conjecture that it is 2-strongly normal."

However, when comparing these pictures with the above Fig. \ref{Oscilat1}, the following conjecture seems more plausible:

The seeming randomness of the digits of $\pi$, $e$, $\sqrt{2}$ and other ``analytically" defined irrational numbers is only a limited one.

\bibliographystyle{splncs03}

\phantomsection
\addcontentsline{toc}{chapter}{References}

\end{document}